\documentclass[12pt]{article}
\usepackage{amsmath}
\usepackage[usenames]{color}
\usepackage{mathrsfs}
\usepackage{amsfonts}
\usepackage{amssymb,amsmath}
\usepackage{CJK}
\usepackage{cite}
\usepackage{cases}
\usepackage{amsthm}

\def\d{\delta}
\def\cl{\centerline}

\def\al{\alpha}
\def\vs{\vspace*}
\def\SV{\mathcal{SV}}
\def\Z{\mathbb{Z}}
\def\C{\mathbb{C}}

\def\QED{\hfill$\Box$}

\numberwithin{equation}{section}
\newtheorem{theo}{Theorem}[section]
\newtheorem{defi}[theo]{Definition}

\newtheorem{lemm}[theo]{Lemma}
\newtheorem{prop}[theo]{Proposition}
\def\G{\Gamma}
\begin{document}
\begin{CJK*}{GBK}{song}

\begin{center}
{\bf\large Structures of not-finitely graded Lie superalgebras\\ related to generalized super-Virasoro algebras\,$^*$}
\footnote {$^{\,*}$Supported by the National Natural Science Foundation of China (No.~11431010, 11371278).

$^{\,\ddag}$ Corresponding author: G. Fan.
}
\end{center}

\cl{Juanjuan Li$^{\,\S,\,\dag}$, Guangzhe Fan$^{\,\S,\,\ddag}$
}

\cl{\small $^{\,\dag}$lijuan103@126.com}
\cl{\small $^{\,\ddag}$yzfanguangzhe@126.com}
\cl{\small $^{\,\S}$Department of Mathematics, Tongji University, Shanghai 200092,
P. R. China}
\vs{8pt}

{\small
\parskip .005 truein
\baselineskip 3pt \lineskip 3pt

\noindent{{\bf Abstract:}  This paper is devoted to investigating the structure theory of a class of not-finitely graded Lie superalgebras related to generalized super-Virasoro algebras. In particular, we completely determine the derivation algebras, the automorphism groups and the second cohomology groups
 of these Lie superalgebras. \vs{5pt}

}

}
\parskip .001 truein\baselineskip 6pt \lineskip 6pt

\section{Introduction}
\setcounter{equation}{0}
For the readers' convenience, we give some notations used in this paper. Let $\Z\subseteq\Gamma$ be an additive subgroup of $\C$ and $s\in \C$  such that $2s\in\Gamma$. For simplicity, let $\Omega$ be an additive subgroup of $\C$ generated by $\Gamma\cup \{s\}$. Denote by $\C$, $\C^*$, $\Z$, $\Z_+$, $\G^*$, $\Omega^*$ the sets of complex numbers, nonzero complex numbers, integers, nonnegative integers, nonzero elements of $\G$, nonzero elements of $\Omega$, respectively. We assume that all vector spaces are based on $\C$, unless otherwise stated. For convenience, the degree of $x$ or $\phi$ is denoted by $|x|$ or $|\phi|$. In addition, $x$ is always assumed to be homogeneous when $|x|$ occurs. Let $L$ be a Lie superalgebra, and denote by $hg(L)$ the set of all homogeneous elements of $L$.

Lie superalgebras as generalizations of Lie algebras were originated from supersymmetry in mathematical physics. The theory of Lie superalgebras plays a prominent role in modern mathematics and physics. In recent years, structures of all kinds of Lie superalgebras have aroused many scholars' great interests (see, e.g., \cite{DHS,KY,H,WPD}). In this paper, we shall investigate the structure theory of a class of not-finitely graded Lie superalgebras related to the generalized super-Virasoro algebras (namely, derivations, automorphisms, 2-cocycles).

Recently, some researchers have studied structures of some Lie algebras related to the Virasoro algebra (see, e.g., \cite{CHS,DHS,DZ,PZ,SXZ,SZ,WWY}). The (generalized) super-Virasoro algebra is closely related to the conformal field theory and the string theory. It plays a very important role in mathematics and physics. The structure and representation theories of (generalized) super-Virasoro algebra have been extensively undertaken by many authors (see, e.g., \cite{DHS,KY,S,S1,S2,S3,SZ1,Y,Z}). The generalized super-Virasoro algebra $SVir[\Gamma,s]$ is a Lie superalgebra whose even part $SV_{\bar{0}}$ has a basis $\{L_{\alpha},C~|~\alpha\in \Gamma\}$ and odd part $SV_{\bar{1}}$ has a basis $\{G_{\mu}~|~\mu\in \Gamma+s\}$, equipped with the following Lie super-brackets:
\begin{eqnarray}\label{1000}
&\!\!\!\!\!\!\!\!\!\!\!\!\!\!\!\!\!\!&
[L_{\alpha},L_{\beta}]=(\beta-\alpha)L_{\alpha+\beta}+\frac{\alpha^3-\alpha}{12}\delta_{\alpha+\beta,0}C, \\
\label{1001}
&\!\!\!\!\!\!\!\!\!\!\!\!\!\!\!\!\!\!&
[L_{\alpha},G_{\mu}]=(\mu-\frac{\alpha}{2})G_{\alpha+\mu},\\
\label{1002}&\!\!\!\!\!\!\!\!\!\!\!\!\!\!\!\!\!\!&
[G_{\mu},G_{\nu}]=2L_{\mu+\nu}+\frac{1}{3}(\mu^2-\frac{1}{4})\delta_{\mu+\nu,0}C,\\
\label{1003}&\!\!\!\!\!\!\!\!\!\!\!\!\!\!\!\!\!\!&
[SV,C]=0,
\end{eqnarray}
where $\alpha, \beta\in\Gamma, \mu,\nu\in{s+\Gamma}$.
If $\Gamma=\Z$, then $SVir[\Z,s](s=0~or~\frac{1}{2})$ is the super-Virasoro algebra. Obviously, $SVir[\Gamma,s]$ is a generalization of the super-Virasoro algebra. Its theories contain the results of the classical super-Virasoro algebra.

In the present work, we will consider the following Lie superalgebra, called not-finitely graded {\it generalized super-Virasoro algebra $\SV[\Gamma,s]$},
which has a basis $\{L_{\al,i},G_{\mu,j}\,|\,\alpha \in\Gamma,
\mu\in {s+\Gamma},\,i,j\in\Z_+\}$, and satisfies the following relations:
\begin{eqnarray}\label{algebra-rela1}
&\!\!\!\!\!\!\!\!\!\!\!\!\!\!\!\!\!\!&
[L_{\al,i},L_{\beta,j}]=
(\beta-\al)L_{\al+\beta,i+j}+(j-i)L_{\al+\beta,i+j-1}, \\
\label{algebra-rela2}
&\!\!\!\!\!\!\!\!\!\!\!\!\!\!\!\!\!\!&
[L_{\al,i},G_{\mu,j}]=
(\mu-\frac{\al}{2}) G_{\al+\mu,i+j}+(j-\frac{i}{2}) G_{\al+\mu,i+j-1},\\
\label{algebra-rela3}&\!\!\!\!\!\!\!\!\!\!\!\!\!\!\!\!\!\!&
[G_{\mu,i},G_{\nu,j}]=2L_{\mu+\nu,i+j},
\end{eqnarray}
where $\alpha, \beta\in\Gamma, \mu,\nu\in{s+\Gamma}, i,j\in\Z_{+}$. We simply denote $\SV=\SV[\Gamma,s]$ . It is obvious that the center of $\SV$ is trival.
 The Lie algebra $W$ with basis
$\{L_{\al,i}\,|\,\al\in\G,\,i\in\Z_+\}$ and relations
\eqref{algebra-rela1}, is the generalized Witt algebra ${\cal W}=W(0,1,0;\G)$ of Witt type studied in \cite{SZ}. Furthermore, the superalgebra $\SV$ contains the centerless generalized super-Virasoro algebra as its subalgebra.

If an undefined notation appears in an expression, we treat it as zero. For example, $L_{\alpha}=0, G_{\mu}=0$ if $\alpha\notin \Gamma, \mu\notin s+\Gamma$.
The Lie superalgebra $\SV$ is $\Omega$-graded
\begin{equation*}\label{100000}
\SV=\raisebox{-5pt}{${}^{\, \, \displaystyle\oplus}_{\alpha\in \Omega}$}{\SV}_{\alpha},\ \ \
 {\SV}_{\alpha}={\rm span} \{L_{\alpha,i}, G_{\alpha,i}\mid i\in \Z_+\}\mbox{ \ for }\al\in\Omega.
\end{equation*}
Obviously, $\SV=\SV_{\overline{0}} \oplus\SV_{\overline{1}}$ is $\Z_{2}$-graded, where
\begin{equation*}\label{100001}
{\SV}_{\overline{0}}={\rm span} \{L_{\alpha,i}\mid \alpha\in \Omega, i\in \Z_+\},\ \ \
{\SV}_{\overline{1}}={\rm span} \{G_{\mu,i}\mid \mu\in \Omega, i\in \Z_+\}.
\end{equation*}

 Because $\Omega$ may not be finitely generated (as a group), $\SV$ may not be finitely generated as a Lie superalgebra. Hence, some classical techniques (such as those in \cite{F}) cannot be directly applied to this case. We must use
some new techniques in order to deal with problems associated with not-finitely generated Lie superalgebras. In addition, in \cite{CHS,FZY} some authors studied some structures of not-finitely graded Lie algebras. Thus, some methods about those Lie algebras can be applied in this case. However, it seems to us that little has been known about not-finitely graded aspect of Lie superalgebras. It may be useful and meaningful for promoting the development of not-finitely graded Lie superalgebras. This is the main motivation of this article.

The present paper is organized as follows. In Section $2$, we review the basic notions about Lie superalgebras. In Section $3$ and Section $4$, we determine the derivations and automorphism groups of $\SV$, respectively. Finally, the second cohomology groups of $\SV$ are obtained in Section $5$. The main results of this paper are summarized in Theorems 3.6, 4.2 and 5.3.

\section{Preliminaries}
In this section, we shall summarize some basic concepts about Lie superalgebras in \cite{K,S}.

\begin{defi}\label{2000}\rm
A Lie superalgebra is a superalgebra $L=L_{\bar{0}}\oplus L_{\bar{1}}$ with multiplication $[\cdot,\cdot]$ satisfying the following two axioms:
\begin{equation}\label{200}
\aligned
&skew~~super-symmetry:~~~~[x,y]=-(-1)^{|x||y|}[y,x],\\
&super~~Jacobi~~identity:~~~~[x,[y,z]]=[[x,y],z]+(-1)^{|x||y|}[y,[x,z]],
\endaligned
\end{equation}
for any $x,y\in hg(L)$, $z\in L$.
\end{defi}

Let $L$ be a Lie superalgebra, then the space $gc(L)$ consisting of all the linear transformations on $L$ has a natural $\Z_{2}$-gradation: $gc(L)=gc(L)_{\overline{0}}\oplus gc(L)_{\overline{1}}$, where
$gc(L)_{\alpha}=\{f\in gc(L)~|~f(L_{\beta})\subseteq L_{\alpha+\beta}$, for any $\beta\in\Z_{2}\}$, for any $\alpha\in\Z_{2}$. Now we give the definition of derivations.
\begin{defi}\label{2001}\rm
A derivation of homogenous $\alpha\in\Z_2$ of  $L$  is an $\C$-linear
transformation $D\in gc(L)_{\alpha}$ such that $$D([x,y])=[D(x),y]+(-1)^{|D||x|}[x, D(y)],$$ for any $x,y\in hg(L)$.
\end{defi}

Denote by ${\rm Der}_{\alpha}L$ the set of all derivations of homogenous $\alpha$ of $L$. Then denote by ${\rm Der\,}{L}={\rm Der\,}_{\bar{0}}{L}\oplus{\rm Der\,}_{\bar{1}}{L}$ the derivation algebra of $L$.
A derivation $D\in {\rm Der}(L)$ is called {\em inner}
if there exists an element $z\in L$ such that $D={\rm ad}_z$,
where ${\rm ad}_z$ is a linear map of $L$ sending $y$ to $[z,y]$ for any $y\in L$.
Denote by ${\rm ad\,}{L}$ the set of all inner derivations of $L$.

\begin{defi}\label{2002}\rm
A bilinear form $\psi:L\times L \rightarrow
\C$ is called a  2-{\it cocycle} on $L$ if $\psi$ satisfies
the following conditions:
\begin{eqnarray*}
\aligned
&skew~~super-symmetry:~~~~\psi(x,y)=-(-1)^{|x||y|}\psi(y,x),\\
&super~~Jacobi~~identity:~~~~\psi(x,[y,z])=\psi([x,y],z)+(-1)^{|x||y|}\psi(y,[x,z]),
\endaligned
\end{eqnarray*}
for any $x,y\in hg(L)$, $z\in L$.
\end{defi}

\begin{defi}\label{2003}\rm
For any linear map $f: L \rightarrow
\C$, we define a 2-cocycle $\psi_{f}$ in the following way:
\begin{eqnarray*}
&&\psi_{f}(x,y)=f([x,y]),
\end{eqnarray*} for any $x,y\in L$.
We call it a  2-{\it coboundary} of $L$.
\end{defi}

Denote by $C^{2}(L,\C)$, $B^{2}(L,\C)$ the vector spaces of all 2-cocycles, 2-coboundaries of $L$ respectively.
The quotient space $H^{2}(L,\C)=C^{2}(L,\C)/B^{2}(L,\C)$
is called  the 2-{\it cohomology group} of $L$.

\section{ Derivations of $\SV$}

In this section, we shall determine the derivations of $\SV$.

Denote by ${\rm Hom}_\Z(\Omega,\C)$ the space of group homomorphisms from $\Omega$ to $\C$. For each $\phi\in {\rm Hom}_\Z(\Omega,\C)$, we define $(c\phi)(\gamma)=c\phi(\gamma)$ for any $r\in\Omega$, $c\in\C$. Furthermore, we define a derivation $D_{\phi}$ as follows:
\begin{equation}\label{3000}
D_{\phi}( L_{\al,i})=\phi(\alpha)L_{\al,i},\ \ D_{\phi}( G_{\al,i})=\phi(\alpha)G_{\al,i}
\ \ \ \mbox{for}\ \ \al\in\Omega,\ i\in\Z_+.
\end{equation}
Denote by ${\rm Hom}_\Z(\Omega,\C)$ the corresponding subspace of ${\rm Der\,} \SV$. In particular, since $\phi_0:\al\mapsto\al$ is in ${\rm Hom}_\Z(\Omega,\C)$, we obtain this special derivation
\begin{equation}\label{3001}
D_0=D_{\phi_0}:L_{\al,i}\mapsto \al L_{\al,i},\ \ G_{\al,i}\mapsto \al G_{\al,i} \mbox{ \ for \ }\al\in\Omega,\ i\in\Z_+.
\end{equation}

\begin{lemm}\label{lemm302}\rm
For every $D\in \SV$, we assume
 \begin{equation}\label{3003}
D=\mbox{$\sum\limits_{\gamma\in\Omega}$} D_\gamma,\ \ D_{\gamma}\in({\rm Der\,}\SV)_\gamma \end{equation}
 such that only finitely many $D_\gamma(x)\neq0$ for every $x\in \SV$ $($such a sum in \eqref{3003} is called summable$)$.
\end{lemm}
\noindent{\it Proof~}~For every $D\in{\rm Der\,}\SV$, we suppose $D(x_\alpha)=\mbox{$\sum\limits
_{\beta\in\Omega}$} y_\beta$, where $x_\alpha\in \SV_\alpha$. Define $D_\gamma(x_\alpha)=y_{\alpha+\gamma}$. Thus, $D_\gamma$ is a derivation of $\SV$.\QED

\begin{lemm}\label{lemm303}\rm
For any $D\in {\rm Der\,}\SV$, replacing $D$ with $D-{\rm ad}_y$ for some $y\in \SV$, we get $D(L_{0,0})=0$.
\end{lemm}

\noindent{\it Proof~}~Assume $D(L_{0,0})\!=\!\sum_{\al,j}(a^{1}_{\al,j}L_{\al,j}+a^{2}_{\al,j}G_{\al,j})\in \SV$ for $a^{\Delta}_{\al,j}\in\C$, where $\Delta=1,2$.
For any $\al\in\Omega$, we define $b^{1}_{\al,j},\ b^{2}_{\al,j} \in\C$ inductively on $j\geq0$ by
\begin{equation}\label{3005}
b^{\Delta}_{\alpha,j}=\left\{\begin{array}{llll}
j^{-1}(-a^{\Delta}_{\al,j-1}-\alpha b^{\Delta}_{\al,j-1}),&\mbox{if \ }j\geq1,\\[4pt]0,&\mbox{if \ }j=0,
\end{array}\right.\end{equation}
where $\Delta=1,2$.

Choosing $y=\sum_{\al,j}(b^{1}_{\al,j}L_{\al,j}+b^{2}_{\al,j}G_{\al,j})\in\SV$, it follows that
\begin{equation}\label{3006}
D(L_{0,0})-{\rm ad}_y(L_{0,0})=0.
\end{equation}\QED

\begin{lemm}\label{lemm304}\rm
If $\gamma\in\Omega^*$, $D\in({\rm Der\,}\SV)_\gamma$ and $D(L_{0,0})=0$, then
$D=0$.
\end{lemm}

\noindent{\it Proof~}~Applying $D$ to $[L_{0,0}, L_{\alpha,0}]=\alpha L_{\alpha,0}$ with $\alpha\in\Omega^*$, then we get
 \begin{equation}\label{3008}
[L_{0,0}, D(L_{\alpha,0})]=\alpha D(L_{\alpha,0}).
\end{equation}
Suppose $D(L_{\alpha,0})=\sum_{j\in\Z_+}(c^{1}_j L_{\alpha+\gamma,j}+c^{2}_j G_{\alpha+\gamma,j})$ for $c^{1}_j,c^{2}_j\in\C$, then we obtain
\begin{equation}\label{3009}
\gamma c^{\Delta}_j=-(j+1)c^{\Delta}_{j+1},
\end{equation}for any $j\geq0$, where $\Delta=1,2$. Thus, $D(L_{\alpha,0})=0$ for any $\alpha\in\Omega$.

Applying $D$ to $[L_{0,0}, L_{0,1}]= L_{0,0}$, then we obtain
\begin{equation}\label{3010}
[L_{0,0}, D(L_{0,1})]=D(L_{0,0})=0.
\end{equation}
Suppose that $D(L_{0,1})=\sum_{j\in\Z_+}(A_j^1L_{r,j}+A_j^2G_{r,j})$, where $A_j^1,A_j^2\in \C$,
then $rA_j^\Delta=-(j+1)A_{j+1}^\Delta$ for any $j\geq0$, where $\Delta=1,2$.
Hence, $ D(L_{0,1})=0$. Now applying $D$ to $[L_{0,0}, G_{\alpha,0}]=\alpha G_{\alpha,0}$, then it shows
 \begin{equation}\label{3011}
[L_{0,0}, D(G_{\alpha,0})]=\alpha D(G_{\alpha,0}).
\end{equation}
Similar to compute $D(L_{\alpha,0})$, we deduce that $D(G_{\alpha,0})=0$ for any $\alpha\in\Omega$.

Since $\SV$ can be generated by $\{L_{\alpha,0},L_{0,1},G_{\alpha,0}\mid\alpha\in\Omega\}$, it forces at once that $D=0$.\QED

\begin{lemm}\label{lemm304}\rm
Suppose $D\in({\rm Der_{\overline{0}}\,}\SV)_0$, $D(L_{0,0})=0$, then
$D(L_{\alpha,0})=d_\alpha L_{\alpha,0},D(G_{\alpha,0})=d_\alpha G_{\alpha,0}$, where $d_\alpha\in{\rm Hom}_\Z(\Omega,\C),\alpha\in\Omega $.
\end{lemm}

 \noindent{\it Proof~}~ Now for any
$D\in ({\rm Der_{\overline{0}}\,}\SV)_0$, we assume
 $D(L_{\alpha,0})=\sum_{\alpha,j}d_{\alpha,j} L_{\alpha,j}$, $D(G_{\alpha,0})=\sum_{\alpha,j}e_{\alpha,j} G_{\alpha,j}$, where $d_{\alpha,j},e_{\alpha,j}\in\C$. Applying $D$ to $[L_{0,0}, L_{\alpha,0}]=\alpha L_{\alpha,0}$ and $[L_{0,0}, G_{\alpha,0}]=\alpha G_{\alpha,0}$, then we have
$d_{\alpha,j}=e_{\alpha,j}=0$ for any $j\geq1$.
Simply denote $d_{\alpha,0}$, $e_{\alpha,0}$ by $d_\alpha$, $e_{\alpha}$, respectively. Hence, $D(L_{\alpha,0})=d_\alpha L_{\alpha,0}, D(G_{\alpha,0})=e_\alpha G_{\alpha,0}$ for any $\alpha\in\Omega$. By applying $D$ to $[L_{\alpha,0}, L_{\beta,0}]=(\beta-\alpha)L_{\alpha+\beta,0}$, it has
\begin{equation}\label{3021}
d_\alpha+d_\beta=d_{\alpha+\beta}  \mbox{ \ for \ }\alpha\neq\beta.
\end{equation}

Furthermore, we obtain that the map $\Phi:\alpha\longmapsto d_{\alpha}$ is an element in ${\rm Hom}_\Z(\Omega,\C)$. Therefore, $d_\alpha\in{\rm Hom}_\Z(\Omega,\C)$.

By using the following three equations
\begin{equation}\label{3022}
\left\{\begin{array}{llll}[L_{\alpha,0},G_{\beta,0}]=(\beta-\frac{\alpha}{2})G_{\alpha+\beta,0},\\[4pt]
[G_{\alpha,0},L_{\beta,0}]=(\frac{\beta}{2}-\alpha)G_{\alpha+\beta,0} ,\\[4pt]
[G_{\alpha,0},G_{\beta,0}]=2L_{\alpha+\beta,0},
\end{array}\right.
\end{equation}
we can conclude that
\begin{equation}\label{3023}
\left\{\begin{array}{llll}d_\alpha+e_\beta=e_{\alpha+\beta}   \mbox{ \ for\ }\alpha\neq2\beta, \\[4pt]
e_\alpha+d_\beta=e_{\alpha+\beta}   \mbox{ \ for\ }\beta\neq2\alpha,\\[4pt]
e_\alpha+e_\beta=d_{\alpha+\beta}.
\end{array}\right.
\end{equation}
Then, we get $e_{\alpha}=d_{\alpha}$ for any $\alpha\in\Omega$.

Thus, $D(L_{\alpha,0})=d_\alpha L_{\alpha,0}, D(G_{\alpha,0})=d_\alpha G_{\alpha,0}$ for any $\alpha\in\Omega$. \QED

\begin{lemm}\label{lemm2}\rm  If $D(L_{0,0})=0$, then
\begin{equation}\label{lemma305}
({\rm Der\,}\SV)_0=({\rm ad\,}\SV)_0\oplus{\rm Hom}_\Z(\Omega,\C).
\end{equation}
\end{lemm}

 \noindent{\it Proof~}~For any $D\in ({\rm Der_{\overline{0}}\,}\SV)_0$, replacing $D'$ by $D-D_\phi$ (cf. (3.1)), then we assume $d_\alpha=0$ for any $\alpha\in\Omega$, thus we get
$D'(L_{\alpha,0})=D'(G_{\alpha,0})=0$.

Now assume $D'(L_{\alpha,i})=\sum_{j\in\Z_+}f_{\alpha,j}^iL_{\alpha,j}$, where $f_{\alpha,j}^i\in \C, i\in \Z_+$. Applying $D'$ to $[L_{0,0},L_{\alpha,1}]=\alpha L_{\alpha,1}+L_{\alpha,0}$, we see $f_{\alpha,j}^1=0$ for any $j\geq1$. Thus, $D'(L_{\alpha,1})=f_{\alpha,0}^1L_{\alpha,0}$ for some $f_{\alpha,0}^1\in \C$. Similarly, applying $D'$ to $[L_{0,0},L_{\alpha,2}]=\alpha L_{\alpha,2}+2L_{\alpha,1}$, we obtain
$$D'(L_{\alpha,2})=f_{\alpha,0}^2L_{\alpha,0}+2f_{\alpha,0}^1L_{\alpha,1}.$$

In view of $[L_{0,1},L_{\alpha,1}]=\alpha L_{\alpha,2}$, it follows
\begin{equation}
\left\{\begin{array}{llll}\alpha f_{0,0}^1=\alpha f_{\alpha,0}^1,\\[4pt]
f_{0,0}^1-f_{\alpha,0}^1=\alpha f_{\alpha,0}^2.
\end{array}\right.
\end{equation}
Thus, $f_{\alpha,0}^1=f_{0,0}^1, f_{\alpha,0}^2=0$ for $\alpha\neq0$. Furthermore, according to $[L_{-\alpha,1},L_{\alpha,1}]=2{\alpha}L_{0,2}$, we get $f_{0,0}^2=0.$ Hence,
$$D'(L_{\alpha,1})=f_{0,0}^1L_{\alpha,0},  D'(L_{\alpha,2})=2f_{0,0}^1L_{\alpha,1}.$$

Induction on $i$, and use the relations of $W$, then we obtain
$$
D'(L_{\alpha,i})=if_{0,0}^1L_{\alpha,i-1} \mbox{ \ for \ }\alpha\in\Omega, i\in \Z_+.
$$
Applying $D'$ to $[G_{0,0},G_{\alpha,i}]=2L_{\alpha,i}$, we have
$$D'(G_{\alpha,i})=if_{0,0}^1G_{\alpha,i} \mbox{ \ for \ }\alpha\in\Omega, i\in \Z_+ .$$

If we replace $y$ (cf. (3.5)) by $y+f_{0,0}^1L_{0,0}$, denote $D'+f_{0,0}^1D_0$ (cf. (3.2)), then $D'(L_{\alpha,i})=D'(G_{\alpha,i})=0$ for $\alpha\in\Omega, i\in \Z_+$.\\

Obviously, $({\rm Der\,}\SV)_0=({\rm Der_{\overline{0}}\,} \SV)_0\oplus ({\rm Der_{\overline{1}}\,}\SV)_0$.

Next for any $D\in({\rm Der_{\overline{1}}\,}\SV)_0$, we suppose
$D(L_{\alpha,i})=\sum_{j\in \Z_+}x_{\alpha,j}^iG_{\alpha,j}, D(G_{\alpha,i})=\sum_{j\in \Z_+}y_{\alpha,j}^iL_{\alpha,j}$,
where $x_{\alpha,j}^i, y_{\alpha,j}^i\in\C$.
Using the following two equalities, $[L_{0,0},L_{\alpha,0}]=\alpha L_{\alpha,0}$ and $[L_{0,0},G_{\alpha,0}]=\alpha G_{\alpha,0}$, then we get
$D(L_{\alpha,0})=x_\alpha G_{\alpha,0}, D(G_{\alpha,0})=y_\alpha L_{\alpha,0}$, where $x_\alpha=x_{\alpha,0}^0, y_\alpha=y_{\alpha,0}^0$. Furthermore, using the relations of $\SV$, we obtain
$$
D(L_{\alpha,i})=D(G_{\alpha,i})=0 \mbox{ \ for any \ }\alpha\in\Omega, i\in \Z_+.
$$
\QED

The following theorem can be concluded by Lemmas~3.1-3.5 immediately, which is the main result of this section.
\begin{theo}\label{theo301}\rm
\begin{equation}\label{3002}
{\rm Der\,}\SV =\raisebox{-5pt}{${}^{\, \, \displaystyle\oplus}_{\gamma\in \Omega}$}({{\rm Der\,}\SV})_{\gamma}={\rm ad\,}\SV\oplus{\rm Hom}_\Z(\Omega,\C), \end{equation}
where $({\rm Der\,}\SV)_\gamma\subset{\rm ad\,}\SV$, if $\gamma\in\Omega^*$, and $({\rm Der\,}\SV)_0=({\rm ad\,}\SV)_0\oplus{\rm Hom}_\Z(\Omega,\C)$.
\end{theo}

\section{Automorphism groups of $\SV$}
The aim of this section is to characterize the automorphism groups of $\SV$.

Denote by ${\rm
Aut\,}{W}$, ${\rm
Aut\,}{\SV}$
the automorphism group of $W$, $\SV$, respectively. Let $\chi(\Omega)$ be the set of characters of $\Omega$, i.e., the set of group homomorphisms $\tau:\Omega\to\C^*$. Set $\Omega^{\C^*}=\{c\in\C^*\,|\,c\Omega=\Omega\}$. A group structure on $\chi(\Omega)\times\Omega^{\C^*}$ can be defined by
\begin{equation}\label{group-sss}
(\tau_1,c_1)\cdot(\tau_2,c_2)=(\tau,c_1c_2),\mbox{ \ where \ }\tau:\al\mapsto\tau_1(c_2\al)\tau_2(\al)\mbox{ \ for\ }\al\in\Omega.
\end{equation}Clearly, we kown that $\chi(\Omega)\times\Omega^{\C^*}$ is just the semidirect product $\chi(\Omega)\rtimes\Omega^{\C^*}$ under the action given by $(c\tau)(\alpha)=\tau(c\alpha)$ for all $c\in\Omega^{\C^*}$, $\tau\in\chi(\Omega),$  $\alpha\in\Omega$.

Define a group homomorphism $\phi:(\tau,c)\mapsto\phi_{\tau,c}$ from $\chi(\Gamma)\times\Gamma^{\C^*}$ to ${\rm Aut\,}W$ such that
 $\phi_{\tau,c}$ is the automorphism of $W$ by
 \begin{equation*}\label{Aususus}
\phi_{\tau,c}:L_{\al,i}\mapsto\tau(\al)c^{i-1}L_{c\al,i},
 \end{equation*}
where $\al\in\Gamma,\,i\in\Z_+$.
We also define a group homomorphism $\phi:{\tau,c}\mapsto\phi_{\tau,c}$ from $\chi(\Omega)\times\Omega^{\C^*}$ to ${\rm Aut\,}\SV$ such that
 $\phi_{\tau,c}$ is the automorphism of $\SV$ by
 \begin{equation*}\label{Aususus}
 \phi_{\tau,c}:L_{\al,i}\mapsto\tau(\al)c^{i-1}L_{c\al,i},
 \end{equation*}
\begin{equation*}\label{Aususus}
\phi_{\tau,c}:G_{\al,i}\mapsto\tau(\al)c^{i-\frac1{2}}G_{c\al,i},
 \end{equation*}
where $\al\in\Omega,\,i\in\Z_+$.

According to \cite{FZY}, we get the following lemma.
\begin{lemm}\label{lemm4001}
${\rm Aut\,}W\cong\chi(\Gamma)\times\Gamma^{\C^*}$.
\end{lemm}

\begin{theo}\label{theo401}\rm
${\rm Aut\,}\SV\cong\chi(\Omega)\times\Omega^{\C^*}\times\Z_{2}$.
\end{theo}

\noindent{\it Proof~}~Assume $\sigma\in {\rm Aut\,}\SV$, $x\in \SV_{\bar{0}}$, we assume $\sigma(x)=y+z$, where $y\in\SV_{\bar{0}}$, $z\in\SV_{\bar{1}}$.
$\sigma([x,x])=[\sigma(x),\sigma(x)]=[z,z]=0$, then $z=0$.
Thus, we have $\sigma(L_{\alpha,i})=\tau(\al)c^{i-1}L_{c\al,i}$ for any $\alpha\in\Omega, i\in \Z_+$.

Obviously, $\sigma(L_{0,0})=c^{-1}L_{0,0}$. Assume $\sigma(G_{\alpha,i})=\Sigma_{\beta,j}(A^{\alpha,i}_{\beta,j}L_{\beta,j}+B^{\alpha,i}_{\beta,j}G_{\beta,j})$,
 where $A^{\alpha,i}_{\beta,j},B^{\alpha,i}_{\beta,j}\in\C$.
Applying $\sigma$ to $[L_{0,0}, G_{\alpha,0}]=\alpha G_{\alpha,0}$, then $[c^{-1}L_{0,0}, \Sigma_{\beta,j}(A^{\alpha,0}_{\beta,j}L_{\beta,j}+B^{\alpha,0}_{\beta,j}G_{\beta,j})]=\alpha \Sigma_{\beta,j}(A^{\alpha,0}_{\beta,j}L_{\beta,j}+B^{\alpha,0}_{\beta,j}G_{\beta,j})$. Comparing the coefficients of $L_{\beta,j}$ and $G_{\beta,j}$,
thus, $\beta=c\alpha$, $A^{\alpha,0}_{c\alpha,j+1}=B^{\alpha,0}_{c\alpha,j+1}=0$ for $j\in\Z_+$. $A^{\alpha,0}_{\beta,j}=B^{\alpha,0}_{\beta,j}=0$ for $\beta\neq c\alpha, j\in\Z_+$. Hence, $\sigma(G_{\alpha,0})=A_{\alpha}L_{c\alpha,0}+B_{\alpha}G_{c\alpha,0}$, where $A_{\alpha}=A^{\alpha,0}_{c\alpha,0}$, $B_{\alpha}=B^{\alpha,0}_{c\alpha,0}$.

Next we apply $\sigma$ to $[L_{\alpha,0},G_{0,0}]=-\frac1{2}\alpha G_{\alpha,0}$, then
\begin{equation}
\left\{\begin{array}{llll}\alpha A_\alpha=2\alpha \tau(\alpha)A_{0},\\[4pt]
\alpha B_\alpha=\alpha \tau(\alpha)B_{0}.
\end{array}\right.
\end{equation}
whence, $A_\alpha=2 \tau(\alpha)A_{0}$, for $\alpha\in\Omega^*$ and $B_\alpha=\tau(\alpha)B_{0}$ for $\alpha\in\Omega$.

Using the equality $[G_{0,0},G_{\alpha,0}]=2L_{\alpha,0}$, we have
\begin{equation}
\left\{\begin{array}{llll}c \alpha A^{2}_{0}+B^{2}_{0}=c^{-1},\\[4pt]
\alpha A_{0}B_{0}=0.
\end{array}\right.
\end{equation}
Thus, $A_{0}=0, B_{0}=\nabla c^{-\frac1{2}}$, where $\nabla=1$ or $-1$.

 Hence, $\sigma(G_{\alpha,0})=\nabla c^{-\frac1{2}}\tau(\alpha)G_{c\alpha,0}$ for $\alpha\in\Omega$.

Applying $\sigma$ to $[G_{0,0},G_{\alpha,i}]=2L_{\alpha,i}$, then $\sigma(G_{\alpha,i})=A^{\alpha,i}_{0,0}L_{0,0}+\nabla c^{i-\frac1{2}}\tau(\alpha)G_{c\alpha,i}$ for any $\alpha \in\Omega, i\in \Z_+$. Thanks to $[L_{0,0},G_{\alpha,i}]=\alpha G_{\alpha,i}+i G_{\alpha,i-1}$, then we get $A^{\alpha,i}_{0,0}=0$ for any $\alpha \in\Omega, i\in \Z_+$.

Finally, $\sigma (L_{\alpha,i})=c^{i-1}\tau(\alpha)L_{c\alpha,i}$, $\sigma (G_{\alpha,i})= \nabla c^{i-\frac1{2}}\tau(\alpha )G_{c\alpha,i}$, where $\alpha\in\Omega,  i\in \Z_+$, $\nabla=1$ or $-1$.

 We have completed the proof.\QED

\section{Second cohomology groups of $\SV$}
Let $\psi\in C^2(\SV,\C)$, we define a linear function $f:\SV\to\C$ such that $f(L_{\al,i})$ and $f(G_{\al,i})$ are given inductively on $i$ as follows
\begin{equation}\label{5000}
f(L_{\al,i})=\left\{\begin{array}{lll}
\frac1{i+1}\psi(L_{0,0},L_{0,i+1}),&\mbox{if \ }\alpha=0,\\[4pt]
\frac1\alpha\Big(\psi(L_{0,0},L_{\alpha,i})-if(L_{\alpha,i-1})\Big),&\mbox{if \ }\alpha\neq0.
\end{array}\right.
\end{equation}
\begin{equation}\label{5001}
f(G_{\al,i})=\left\{\begin{array}{lll}
\frac{2}{2i-1}\psi(L_{0,1},G_{0,i}),&\mbox{if \ }\alpha=0,\\[4pt]
\frac1\alpha\Big(\psi(L_{0,0},G_{\alpha,i})-if(G_{\alpha,i-1})\Big),&\mbox{if \ }\alpha\neq0.
\end{array}\right.
\end{equation}Set $\phi=\psi-\psi_f$. According to \cite{SZ}, we obtain the following formula:
\begin{equation}\label{5002}
~~\phi(L_{\al,i},L_{\beta,j})=0,
\end{equation}
for any $\al,\beta\in\Omega,\, i,j\in\Z_+.$

\begin{prop}\label{theo500}\rm We have
\begin{equation}\label{5003}
~~\phi(L_{\al,i},G_{\beta,j})=0,
\end{equation}
for any $\al,\beta\in\Omega,\, i,j\in\Z_+.$
 \end{prop}
\noindent{\it Proof~}~Firstly, we have
\begin{eqnarray}
\label{5004}
\phi(L_{0,0},G_{\beta,j})&\!\!\!=\!\!\!&\psi(L_{0,0},G_{\beta,j})\!-\!f([L_{0,0},G_{\beta,j}])\!=\!0
\mbox{ \ for \ } \beta\neq0.
\end{eqnarray}

One has
\begin{eqnarray*}\!\!\!\!\!\!\!\!\!\!\!\!&\!\!\!\!\!\!\!\!\!\!\!\!\!\!\!&
\ \ \ \ \ \phi(L_{0,0},G_{0,0})=2\phi(L_{0,0},[L_{0,1},G_{0,0}])
\nonumber\\\!\!\!\!\!\!\!\!\!\!\!\!&\!\!\!\!\!\!\!\!\!\!\!\!\!\!\!&
 \ \ \ \ \ \ \ \ \ \ \ \ \ \ \ \ \ \ \ \ \ \,=2(\phi(L_{0,0},G_{0,0})+\phi(L_{0,1},0))
\nonumber\\\!\!\!\!\!\!\!\!\!\!\!\!&\!\!\!\!\!\!\!\!\!\!\!\!\!\!\!&
 \ \ \ \ \ \ \ \ \ \ \ \ \ \ \ \ \ \ \ \ \ \,=2\phi(L_{0,0},G_{0,0}).
\end{eqnarray*}
It is clear that $\phi(L_{0,0},G_{0,0})=0$.

Furthermore, for any $j\geq 1$, we obtain
\begin{eqnarray*}\!\!\!\!\!\!\!\!\!\!\!\!&\!\!\!\!\!\!\!\!\!\!\!\!\!\!\!&
\ \ \ \ \ \ \ \ \phi(L_{0,0},G_{0,j})=\psi(L_{0,0},G_{0,j})-jf(G_{0,j-1})
\nonumber\\\!\!\!\!\!\!\!\!\!\!\!\!&\!\!\!\!\!\!\!\!\!\!\!\!\!\!\!&
\ \ \ \ \ \ \ \ \ \ \ \ \ \ \ \ \ \ \ \ \ \ \ \ \, =\psi(L_{0,0},G_{0,j})-\frac{2j}{2j-3}\psi(L_{0,1},G_{0,j-1})
\nonumber\\\!\!\!\!\!\!\!\!\!\!\!\!&\!\!\!\!\!\!\!\!\!\!\!\!\!\!\!&\ \ \ \ \ \ \ \ \ \ \ \ \ \ \ \ \ \ \ \ \ \ \  \ \,
=\psi(L_{0,0},G_{0,j})-\frac{2}{2j-3}\psi(L_{0,1},[L_{0,0},G_{0,j}])
\nonumber\\\!\!\!\!\!\!\!\!\!\!\!\!&\!\!\!\!\!\!\!\!\!\!\!\!\!\!\!&\ \ \  \ \ \ \ \ \ \ \ \ \ \ \ \ \ \ \ \ \ \ \ \ \,
=\psi(L_{0,0},G_{0,j})-\frac{2}{2j-3}(-\psi(L_{0,0},G_{0,j})+\frac{2j-1}{2}\psi(L_{0,0},G_{0,j}))
\nonumber\\\!\!\!\!\!\!\!\!\!\!\!\!&\!\!\!\!\!\!\!\!\!\!\!\!\!\!\!&\ \ \ \ \ \ \ \ \ \ \ \ \ \ \ \ \ \ \ \ \ \ \ \ \,=0.
\end{eqnarray*}
Thus,
\begin{equation}
\label{5005}
\phi(L_{0,0},G_{\beta,j})=0 \ \ \ \ \ \ \mbox{ for } \ \beta\in\Omega,\ j\in\Z_+.
\end{equation}

For $\alpha+\beta\neq0$, it shows
\begin{eqnarray*}\!\!\!\!\!\!\!\!\!\!\!\!&\!\!\!\!\!\!\!\!\!\!\!\!\!\!\!&
\ \ \ \ \ \ \ \ \phi(L_{\alpha,0},G_{\beta,0})=\psi(L_{\alpha,0},G_{\beta,0})-(\beta-\frac{\alpha}{2}) f(G_{\alpha+\beta,0})
\nonumber\\\!\!\!\!\!\!\!\!\!\!\!\!&\!\!\!\!\!\!\!\!\!\!\!\!\!\!\!&
\ \ \ \ \ \ \ \ \ \ \ \ \ \ \ \ \ \ \ \ \ \ \ \ \  =\psi(L_{\alpha,0},G_{\beta,0})- \frac{2\beta-\alpha}{2(\alpha+\beta)}\psi(L_{0,0},G_{\alpha+\beta,0})
\nonumber\\\!\!\!\!\!\!\!\!\!\!\!\!&\!\!\!\!\!\!\!\!\!\!\!\!\!\!\!&
\ \ \ \ \ \ \ \ \ \ \ \ \ \ \ \ \ \ \ \ \ \ \ \ \  =\psi(L_{\alpha,0},G_{\beta,0})-\frac{1}{\alpha+\beta}\psi(L_{0,0},[L_{\alpha,0},G_{\beta,0}])
\nonumber\\\!\!\!\!\!\!\!\!\!\!\!\!&\!\!\!\!\!\!\!\!\!\!\!\!\!\!\!&
\ \ \ \ \ \ \ \ \ \ \ \ \ \ \ \ \ \ \ \ \ \ \ \ \  =0.
\end{eqnarray*}
By induction on $i+j$ for $\alpha+\beta\neq0$, suppose $\phi(L_{\alpha,i-1},G_{\beta,j})=\phi(L_{\alpha,i},G_{\beta,j-1})=0$.
Applying the Jacobi identity on $(L_{0,0},L_{\alpha,i},G_{\beta,j})$, we have
\begin{eqnarray*}\!\!\!\!\!\!\!\!\!\!\!\!&\!\!\!\!\!\!\!\!\!\!\!\!\!\!\!&
\ \ \ \ \ \ \ \ 0=\phi(L_{0,0},[L_{\alpha,i},G_{\beta,j}])
\nonumber\\\!\!\!\!\!\!\!\!\!\!\!\!&\!\!\!\!\!\!\!\!\!\!\!\!\!\!\!&
 \ \ \ \ \ \ \ \ \ \ \,=\phi([L_{0,0},L_{\alpha,i}],G_{\beta,j})+\phi(L_{\alpha,i},[L_{0,0},G_{\beta,j}])
\nonumber\\\!\!\!\!\!\!\!\!\!\!\!\!&\!\!\!\!\!\!\!\!\!\!\!\!\!\!\!&
 \ \ \ \ \ \ \ \ \ \ \,=(\alpha+\beta)\phi(L_{\alpha,i},G_{\beta,j}).
\end{eqnarray*}
Hence, we suppose
\begin{equation}
\label{510000}
\phi(L_{\al,i},G_{\beta,j})=\d_{\al+\beta,0}\phi(L_{\alpha,i},G_{-\alpha,j})\ \ \ \ \ \mbox{ for } \alpha,\beta\in\Omega,\ i,j\in\Z_+.
\end{equation}

Next we shall compute $\phi(L_{\al,i},G_{-\al,j})$. It follows that
\begin{eqnarray*}\!\!\!\!\!\!\!\!\!\!\!\!&\!\!\!\!\!\!\!\!\!\!\!\!\!\!\!&
\ \ \ \ \ \ \ \ \ \  \ \ \ \ \ \ \ \ \ \,0=\phi(L_{0,0},[L_{\al,0},G_{-\al,j}])
\nonumber\\\!\!\!\!\!\!\!\!\!\!\!\!&\!\!\!\!\!\!\!\!\!\!\!\!\!\!\!&
 \ \ \ \ \ \ \ \ \ \ \ \ \ \ \ \ \ \ \ \ \ \,=\phi([L_{0,0},L_{\alpha,0}],G_{-\alpha,j})+\phi([G_{-\alpha,j},L_{0,0}],L_{\alpha,0})
\nonumber\\\!\!\!\!\!\!\!\!\!\!\!\!&\!\!\!\!\!\!\!\!\!\!\!\!\!\!\!&
 \ \ \ \ \ \ \ \ \ \ \ \ \ \ \ \ \ \ \ \ \ \,=\al\phi(L_{\al,0},G_{-\al,j})-\al\phi(L_{\al,0},G_{-\al,j})+j\phi(L_{\al,0},G_{-\al,j-1})
 \nonumber\\\!\!\!\!\!\!\!\!\!\!\!\!&\!\!\!\!\!\!\!\!\!\!\!\!\!\!\!&
 \ \ \ \ \ \ \ \ \ \ \ \ \ \ \ \ \ \ \ \ \ \,=j\phi(L_{\al,0},G_{-\al,j-1}).
\end{eqnarray*}
Thus, $\phi(L_{\al,0},G_{-\al,j})=0$ for $\al\in\Omega$, $j\in\Z_{+}$.

Induction on $i$, suppose $\phi(L_{\al,i-1},G_{-\al,j})=0$, then
\begin{eqnarray*}\!\!\!\!\!\!\!\!\!\!\!\!&\!\!\!\!\!\!\!\!\!\!\!\!\!\!&
\ \ \ \ \ \ \ \ \ \  \ \ \ \ \ \ \ \ \ \,0=\phi(L_{0,0},[L_{\al,i},G_{-\al,j}])
\nonumber\\\!\!\!\!\!\!\!\!\!\!\!\!&\!\!\!\!\!\!\!\!\!\!\!\!\!\!\!&
 \ \ \ \ \ \ \ \ \ \ \ \ \ \ \ \ \ \ \ \ \ \,=\phi([L_{0,0},L_{\al,i}],G_{-\al,j})+\phi([G_{-\alpha,j},L_{0,0}],L_{\alpha,i})
\nonumber\\\!\!\!\!\!\!\!\!\!\!\!\!&\!\!\!\!\!\!\!\!\!\!\!\!\!\!\!&
 \ \ \ \ \ \ \ \ \ \ \ \ \ \ \ \ \ \ \ \ \ \,=\alpha\phi(L_{\alpha,i},G_{-\al,i})+i\phi(L_{\alpha,i-1},G_{-\alpha,j})-\alpha\phi(L_{\alpha,i},G_{-\al,j})+j\phi(L_{\alpha,i},G_{-\alpha,j-1})
\nonumber\\\!\!\!\!\!\!\!\!\!\!\!\!&\!\!\!\!\!\!\!\!\!\!\!\!\!\!\!&
 \ \ \ \ \ \ \ \ \ \ \ \ \ \ \ \ \ \ \ \ \ \,=j\phi(L_{\alpha,i},G_{-\al,j-1}).
\end{eqnarray*}
Thus, $\phi(L_{\al,i},G_{-\al,j})=0$ for any $i,j\in\Z_+$, $\al\in\Omega$.
This implies that $\phi(L_{\al,i},G_{\beta,j})=0$ for $\al,\beta\in\Omega$, $j\in\Z_{+}$.\QED

\begin{prop}\label{theo5.3}\rm We have
\begin{equation}\label{58000}
\phi(G_{\al,i},G_{\beta,j})=0,
\end{equation}
for any $\alpha,\beta\in\Omega$, $i,j\in\Z_+$.
\end{prop}
\noindent{\it Proof~}~According to the equation~(5.3), we deduce that $\phi(L_{0,0},L_{\alpha,i})=0$  for any $\al \in\Omega$, $i\in\Z_{+}$.

We also obtain that
\begin{eqnarray*}\!\!\!\!\!\!\!\!\!\!\!\!&\!\!\!\!\!\!\!\!\!\!\!\!\!\!\!&
\ \ \ \ \ \ \  \ \ \ \ \  0=\phi(L_{0,0},2L_{\alpha+\beta,0})
\nonumber\\\!\!\!\!\!\!\!\!\!\!\!\!&\!\!\!\!\!\!\!\!\!\!\!\!\!\!\!&
 \ \ \ \ \ \ \ \ \ \ \ \ \ \ \,=\phi(L_{0,0},[G_{\alpha,0},G_{\beta,0}])
\nonumber\\\!\!\!\!\!\!\!\!\!\!\!\!&\!\!\!\!\!\!\!\!\!\!\!\!\!\!\!&
\ \ \ \ \ \ \ \ \ \ \ \ \ \ \,=(\alpha+\beta)\phi(G_{\al,0},G_{\beta,0}).
\end{eqnarray*}

It turns out that $\phi(G_{\al,0},G_{\beta,0})=0$ for $\alpha+\beta\neq0$. Induction on $i+j$, suppose $\phi(G_{\al,i-1},G_{\beta,j})=\phi(G_{\al,i},G_{\beta,j-1})=0$. Then we get
\begin{eqnarray*}\!\!\!\!\!\!\!\!\!\!\!\!&\!\!\!\!\!\!\!\!\!\!\!\!\!\!\!&
\ \ \ \ \ \ \  \ \ \ \ \ 0=\phi(L_{0,0},2L_{\alpha+\beta,i+j})
\nonumber\\\!\!\!\!\!\!\!\!\!\!\!\!&\!\!\!\!\!\!\!\!\!\!\!\!\!\!\!&
 \ \ \ \ \ \ \ \ \ \ \ \ \ \ \,=\phi(L_{0,0},[G_{\alpha,i},G_{\beta,j}])
\nonumber\\\!\!\!\!\!\!\!\!\!\!\!\!&\!\!\!\!\!\!\!\!\!\!\!\!\!\!\!&
\ \ \ \ \ \ \ \ \ \ \ \ \ \ \,=(\alpha+\beta)\phi(G_{\al,i},G_{\beta,j})+j\phi(G_{\al,i},G_{\beta,j-1})+i\phi(G_{\al,i-1},G_{\beta,j})
\nonumber\\\!\!\!\!\!\!\!\!\!\!\!\!&\!\!\!\!\!\!\!\!\!\!\!\!\!\!\!&
\ \ \ \ \ \ \ \ \ \ \ \ \ \ \,=(\alpha+\beta)\phi(G_{\alpha,i},G_{\beta,j}).
\end{eqnarray*}
Therefore, $\phi(G_{\al,i},G_{\beta,j})=0$ for $\alpha+ \beta\neq0$. So in the following, we shall only be concerned with the case $\alpha+ \beta=0$.\\

If $j\in\Z_+,$ then
\begin{eqnarray*}\!\!\!\!\!\!\!\!\!\!\!\!&\!\!\!\!\!\!\!\!\!\!\!\!\!\!\!&
\ \ \ \ \ \ \  \phi(G_{0,0},G_{0,j})=\frac1{j+1}\phi(G_{0,0},[L_{0,0},G_{0,j+1}])
\nonumber\\\!\!\!\!\!\!\!\!\!\!\!\!&\!\!\!\!\!\!\!\!\!\!\!\!\!\!\!&
\ \ \ \ \ \ \ \ \ \ \ \ \ \ \ \ \ \ \ \ \ \ \ \,=\frac1{j+1}(\phi([G_{0,0},L_{0,0}],G_{0,j+1})+\phi(L_{0,0},[G_{0,0},G_{0,j+1}]))
\nonumber\\\!\!\!\!\!\!\!\!\!\!\!\!&\!\!\!\!\!\!\!\!\!\!\!\!\!\!\!&
\ \ \ \ \ \ \ \ \ \ \ \ \ \ \ \ \ \ \ \ \ \ \ \,=\frac1{j+1}\phi(L_{0,0},2L_{0,j+1}).
\nonumber\\\!\!\!\!\!\!\!\!\!\!\!\!&\!\!\!\!\!\!\!\!\!\!\!\!\!\!\!&
\ \ \ \ \ \ \ \ \ \ \ \ \ \ \ \ \ \ \ \ \ \ \ \,=0.
\end{eqnarray*}
Induction on $i$, suppose $\phi(G_{0,i-1},G_{0,j})=0$ for any $j\in\Z_+$, then
\begin{eqnarray*}\!\!\!\!\!\!\!\!\!\!\!\!&\!\!\!\!\!\!\!\!\!\!\!\!\!\!\!&
\ \ \ \ \ \ \ \ \ \ \ \ 0=\phi(L_{0,0},2L_{0,i+j})
\nonumber\\\!\!\!\!\!\!\!\!\!\!\!\!&\!\!\!\!\!\!\!\!\!\!\!\!\!\!\!&
 \ \ \ \ \ \ \ \ \ \ \ \ \ \ \,=\phi(L_{0,0},[G_{0,i},G_{0,j}])
\nonumber\\\!\!\!\!\!\!\!\!\!\!\!\!&\!\!\!\!\!\!\!\!\!\!\!\!\!\!\!&
 \ \ \ \ \ \ \ \ \ \ \ \ \ \ \,=\phi([L_{0,0},G_{0,i}],G_{0,j})+\phi(G_{0,i},[L_{0,0},G_{0,j}])
\nonumber\\\!\!\!\!\!\!\!\!\!\!\!\!&\!\!\!\!\!\!\!\!\!\!\!\!\!\!\!&
 \ \ \ \ \ \ \ \ \ \ \ \ \ \ \,=i\phi(G_{0,i-1},G_{0,j})+j\phi(G_{0,i},G_{0,j-1})
\nonumber\\\!\!\!\!\!\!\!\!\!\!\!\!&\!\!\!\!\!\!\!\!\!\!\!\!\!\!\!&
 \ \ \ \ \ \ \ \ \ \ \ \ \ \ \,=j\phi(G_{0,i},G_{0,j-1}).
\end{eqnarray*}
Thus, $\phi(G_{0,i},G_{0,j})=0$ for any $i,j\in\Z_+$. Furthermore, it follows that
\begin{eqnarray*}\!\!\!\!\!\!\!\!\!\!\!\!&\!\!\!\!\!\!\!\!\!\!\!\!\!\!\!&
\ \ \ \ \ \ 0=\phi(L_{0,0},2L_{0,i})
\nonumber\\\!\!\!\!\!\!\!\!\!\!\!\!&\!\!\!\!\!\!\!\!\!\!\!\!\!\!\!&
\ \ \ \ \ \ \ \ \,=\phi(L_{0,0},[G_{\alpha,i},G_{-\alpha,0}])
\nonumber\\\!\!\!\!\!\!\!\!\!\!\!\!&\!\!\!\!\!\!\!\!\!\!\!\!\!\!\!&
\ \ \ \ \ \ \ \ \,=i\phi(G_{\alpha,i-1},G_{-\alpha,0}).
\end{eqnarray*}

Hence, $\phi(G_{\alpha,i},G_{-\alpha,0})=0$ for any $i\in\Z_{+}$.

Applying the Jacobi identity on $(L_{0,0},G_{\alpha,i},G_{-\alpha,j})$, then we get
\begin{eqnarray*}\!\!\!\!\!\!\!\!\!\!\!\!&\!\!\!\!\!\!\!\!\!\!\!\!\!\!\!&
\ \ \ \ \ \ \ \ \ 0=2\phi(L_{0,0},L_{0,i+j})
\nonumber\\\!\!\!\!\!\!\!\!\!\!\!\!&\!\!\!\!\!\!\!\!\!\!\!\!\!\!\!&
\ \ \ \ \ \ \ \ \ \ \ \,=\phi(L_{0,0},[G_{\alpha,i},G_{-\alpha,j}])
\nonumber\\\!\!\!\!\!\!\!\!\!\!\!\!&\!\!\!\!\!\!\!\!\!\!\!\!\!\!\!&
\ \ \ \ \ \ \ \ \ \ \ \,=i\phi(G_{\alpha,i-1},G_{-\alpha,j})+j\phi(G_{\alpha,i},G_{-\alpha,j-1}).
\end{eqnarray*}
Induction on $j$, we also obtain $\phi(G_{\alpha,i},G_{-\alpha,j})=0$.

By previous computations, it will suffice to show the result.\QED

Now we get the main result of this section.

\begin{theo}\label{theo501}\rm The second cohomology group of $\SV$ is trivial, i.e.
\begin{equation}\label{lemma305}
{\rm C^2}(\SV,\C)= {\rm B^2}(\SV,\C).
\end{equation}
\end{theo}

\end{CJK*}
\end{document}